\documentclass{amsart}

\usepackage{amsmath,amssymb,amsfonts}

\newtheorem{theorem}{Theorem}
\newcommand{\bt}{\begin{theorem}}
\newcommand{\et}{\end{theorem}}

\newtheorem{lemma}{Lemma}
\newcommand{\bl}{\begin{lemma}}
\newcommand{\el}{\end{lemma}}

\newtheorem{oproblem}{Open Problem} 
\newcommand{\boprob}{\begin{oproblem}}
\newcommand{\eoprob}{\end{oproblem}}

\newcommand{\beq}{\begin{equation}}
\newcommand{\eeq}{\end{equation}}
\newcommand{\benum}{\begin{enumerate}}
\newcommand{\eenum}{\end{enumerate}}

\newcommand{\N}{\ensuremath{ \mathbf N }}
\newcommand{\Z}{\ensuremath{\mathbf Z}}
\newcommand{\bq}{\begin{eqnarray*}}\newcommand{\eq}{\end{eqnarray*}}\newcommand{\ba}{\begin{array}}\newcommand{\ea}{\end{array}}

\DeclareMathOperator{\card}{\text{card}}

\title{Cassels bases}

\author{Melvyn B. Nathanson}

\thanks{This work was supported in part by the PSC-CUNY Research Award Program.}

\address{Department of Mathematics,
Lehman College (CUNY), Bronx, New York 10468, and
CUNY Graduate Center, New York, New York 10016}

\email{melvyn.nathanson@lehman.cuny.edu}



\begin{document}
\maketitle
\begin{abstract}
This paper describes several classical constructions of thin bases of finite order in additive number theory, and, in particular, gives a complete presentation of a beautiful construction of J. W. S. Cassels of a class of polynomially asymptotic bases.  Some open problems are also discussed.
\end{abstract}

\section{Additive bases of finite order}

The fundamental object in additive number theory is the {\em sumset}.
If $h \geq	2$ and $A_1,\ldots, A_h$ are sets of integers, then we define the sumset
\beq            \label{thin:hsumset}
A_1+\cdots+A_h = \{a_1+\cdots + a_h: a_i \in A_i \mbox{ for } i = 1,\ldots,h\}.
\eeq
If $A_1 = A_2 = \cdots = A_h = A,$ then the sumset
\beq            \label{thin:hA}
hA = \underbrace{A+A+\cdots + A}_{\mbox{$h$ summands}}
\eeq
is called the {\em $h$-fold sumset of $A$.}
If $0 \in A$, then
\[
A \subseteq 2A \subseteq \cdots \subseteq hA \subseteq (h+1)A\subseteq \cdots
\]
For example, 
\[
\{0,1,4,5\} + \{0,2,8,10\} = [0,15]
\]
and
\bq
\lefteqn{\{3,5,7,11\} + \{3,5,7,11,13,17,19\} =}\\
& & \{6,8,10,12,14,16,18,20,22,24,26,28,30\}.
\eq

The set $A$ is called a \emph{basis of order $h$}\index{basis} for the set $B$ if every element of $B$ can be represented as the sum of exactly $h$ not necessarily distinct elements of $A$, or, equivalently, if $B \subseteq hA$.  The set $A$ is an \emph{asymptotic basis of order $h$}\index{asymptotic basis} for $B$ if the sumset $hA$ contains all but finitely many elements of $B$, that is, if $\card(B \setminus hA) < \infty$.   The set $A$ is a basis  (resp. asymptotic basis)
of finite order for $B$ if $A$ is a basis  (resp. asymptotic basis) of order $h$ for $B$ for some positive integer $h$.  The set $A$ of nonnegative integers is a basis of finite order for the nonnegative integers only if $0,1 \in A$.

Many classical results and conjectures in additive number theory state that some ``interesting'' or ``natural'' set of nonnegative integers is a basis or asymptotic basis of finite order.   
 For example, the Goldbach conjecture  \index{Goldbach conjecture} asserts that the set
of odd prime numbers is a basis of order 2 for the even integers
greater than 4.
Lagrange's theorem \index{Lagrange's theorem} states the set of squares is a basis of order 4 for the nonnegative integers $\N_0$.  
Wieferich \index{Wieferich's theorem}proved that  the set of nonnegative cubes is a basis of order 9 for  $\N_0$, and Linnik  \index{Linnik's theorem}proved that the set of nonnegative cubes is an asymptotic  basis of order 7 for  $\N_0$.   More generally, for any integer $k \geq 2,$   \index{Waring's problem}
Waring's problem, proved by Hilbert \index{Hilbert's theorem}in 1909, 
states that  the set of nonnegative $k$-th powers is a basis of finite order for $\N_0$.
Vinogradov  \index{Vinogradov's theorem} proved that the set
of odd prime numbers is an asymptotic basis of order 3 for the odd positive integers.  Nathanson~\cite{nath96aa} contains complete proofs of all of these results.

\emph{Notation:}  
Let \N, $\N_0$, and \Z\ denote the sets of positive integers, nonnegative integers, and integers, respectively.  For real numbers $x$ and $y$, we define the intervals of integers
$[x,y] = \{n\in \Z : x \leq n \leq y\}$,
$(x,y] = \{n\in \Z : x < n \leq y\}$,
and 
$[x,y) = \{n\in \Z : x \leq n < y\}$. 
For any sets $A$ and $A'$ of integers and any integer $c$, we  define the {\em difference set}
\index{difference set}
\[
A - A' = \{ a-a':a\in A \text{ and } a' \in A' \}
\]
and the {\em dilation}\index{dilation}  by $c$ of the set $A$
\[
c\ast A = \{ca:a\in A\}.
\]
Thus, $2\ast\N$ is the set of positive even integers, and $
2\ast \N - \{0,1\} = \N.$

Denote the cardinality of the set $X$ by $|X|$.

Let $f$ be a complex-valued  function on the domain $\Omega$ and let $g$ be a positive function on the domain $\Omega$.  Usually $\Omega$ is the set of positive integers or the set of all real numbers $x \geq x_0$.    We write $f \ll g$ or $f = O(g)$ if there is a number $c> 0$ such that $|f(x)| \leq cg(x)$ for all $x \in \Omega$.  We write $f \gg g$ if  there is a number $c> 0$ such that $|f(x)| \geq cg(x)$ for all $x \in \Omega$. We write $f = o(g)$ if $\lim_{x\rightarrow \infty} f(x)/g(x) = 0$.

\section{A lower bound for bases of finite order}

For any set $A$ of integers, the {\em counting function} of $A$, denoted $A(x)$,
\index{counting function}
counts the number of positive integers in $A$ not exceeding $x$, that is,
\[
A(x) = \sum_{\substack{a \in A \\1 \leq a \leq x}} 1 = \left| A \cap [1,x] \right|.
\]

\bt           \label{thin:theorem:thin1}
Let $h \geq 2$ and let $A = \{a_k\}_{k=1}^{\infty}$ be a set of nonnegative integers with $a_k < a_{k+1}$ for all $k \geq 1$.  
If $A$ is an asymptotic basis of order $h$, then
\beq            \label{thin:thin1}
A(x) \gg x^{1/h}
\eeq
for all sufficiently large real numbers $x$ and
\beq            \label{thin:thin2}
a_k \ll k^h
\eeq
for all positive integers $k$.
If $A$ is a basis of order $h$, then inequality~(\ref{thin:thin1})
 holds for all real numbers $x \geq 1$.
\et

\begin{proof}
If $A$ is an asymptotic basis of order $h$, then there exists an integer $n_0$
such that every integer $m \geq n_0$ can be represented as the sum 
of $h$ elements of $A$.  Let $x \geq x_0$ and let $n$ be the integer part of $x$.   Then $A(x) = A(n)$.  
There are $n-n_0+1$ integers $m$ such that
\[
n_0 \leq m \leq n.
\]
Since the elements of $A$ are nonnegative integers, 
it follows that if
\[
m = a'_1 + \cdots + a'_h \qquad\mbox{with $a'_k \in A$ for $k = 1,\ldots, h$,}
\]
then
\[
0 \leq a'_k \leq m\leq n  \qquad\mbox{for $k = 1,\ldots, h$.}
\]
The set $A$ contains exactly $A(n)$ positive integers not exceeding $n$,
and $A$ might also contain 0, hence $A$ contains 
at most $A(n)+1$ nonnegative integers not exceeding $n$.
Since the number of ways to choose $h$ elements with repetitions 
from a set of cardinality $A(n)+1$ is $\binom{A(n)+h}{h}$, 
it follows that
\[
n+1-n_0 \leq {A(n)+h\choose h} < \frac{(A(n)+h)^h}{h!}
\]
and so
\[
A(x) = A(n) > \left( h!(n+1-n_0)\right)^{1/h}-h \gg n^{1/h}  \gg x^{1/h}
\]
for all sufficiently large  $x$.  
We have $A(a_k) = k$ if $a_1 \geq 1$ and $A(a_k) = k-1$ if $a_1 = 0$, hence 
\[
k \geq  A(a_k) \gg a_k^{1/h}
\]
or, equivalently, 
\[
a_k \ll k^h
\]
for all sufficiently large integers $k$, hence for all positive integers $k$.

If $A$ is a basis of order $h$, then $1 \in A$ and so $A(n)/n > 0$ for all $n \geq 1$.  Therefore, $A(x) \gg x^{1/h}$
for all $x \geq 1.$
This completes the proof.
\end{proof}

Let $A$ be a set of nonnegative integers.
By Theorem~\ref{thin:theorem:thin1}, if $A$ is an asymptotic basis of order $h$,  then $A(x) \gg x^{1/h}.$   If $A$ is an asymptotic basis of order $h$ such that 
\[
A(x) \ll x^{1/h},
\]
then $A$ is called a {\em thin asymptotic basis} of order $h$.\index{thin basis}
\index{basis!thin}
If $hA=\N_0$ and $A(x) \ll x^{1/h},$ then $A$ is called a {\em thin basis} of order $h$.  In the next section we construct examples of thin bases.

\section{Raikov-St{\" o}hr bases}
In 1937 Raikov and St\" ohr independently published the first examples of thin
bases for the natural numbers.
Their construction is based on the fact that every nonnegative integer 
can be written uniquely as the sum of pairwise distinct powers of 2.
The sets constructed in the following theorem 
will be called {\em Raikov-St\" ohr bases}.
\index{Raikov-St\" ohr basis}

\bt[Raikov-St\" ohr]   \label{thin:theorem:Raikov-Stohr}
Let $h \geq 2.$  For $i = 0,1,\ldots, h-1,$ let $W_i = \{ i, h+i, 2h+i,\ldots\}$ 
denote the set of all
nonnegative integers that are congruent to $i$ modulo $h,$
and let $\mathcal{F}(W_i)$ be the set of all finite subsets of $W_i.$
Let
\[
A_i = \left\{ \sum_{f\in F} 2^f : F \in \mathcal{F}(W_i) \right\}
\]
and 
\[
A = A_0 \cup A_1 \cup \cdots \cup A_{h-1}.
\]
Then $A$ is a thin basis of order $h$.
\et

\begin{proof}
Note that for all $i = 0,1,\ldots, h-1$ we have $0 \in A_i$ 
since $\emptyset \in \mathcal{F}(W_i)$ 
and $\sum_{f\in \emptyset} 2^f = 0.$   This implies that 
\[
A_0 + A_1 + \cdots + A_{h-1}  
\subseteq  h\left( \bigcup_{i=0}^{h-1} A_i \right) 
= hA
\]
Moreover, $A_i \cap A_j = \{ 0\}$ if  $0 \leq i < j \leq h-1.$

First we show that $A$ is a basis of order $h$.
Every positive integer $n$ is uniquely the  sum of distinct powers of two, 
so we can write
\[
n = \sum_{j=0}^{\infty}\varepsilon_j 2^j,
\]
where the sequence $\{\varepsilon_j\}_{j=0}^{\infty}$ satisfies
$\varepsilon_j \in \{0,1\}$ for all $j \in \N_0$ and $\varepsilon_j = 0$
for all sufficiently large $j$.
Since
\[
\sum_{ \substack{j=0 \\ j \equiv i\pmod{h} }}^{\infty}\varepsilon_j 2^j \in A_i,
\]
it follows that 
\begin{align*}
n & =  \sum_{j=0}^{\infty}\varepsilon_j 2^j \\
& =  \sum_{i=0}^{h-1} 
\left( \sum_{ \substack{ j=0 \\ j \equiv i\pmod{h}}}^{\infty}\varepsilon_j 2^j \right) \\
& \in  A_0 + A_1 + \cdots + A_{h-1} \\
& \subseteq  hA
\end{align*}
and so $A$ is a basis of order $h$.

We shall compute the counting functions of the sets $A_i$ and $A$.  
Let  $x \geq 2^{h-1}$.  For every $i \in \{ 0,1,\ldots, h-1\}$, there is a unique positive integer $r$ such that 
\[
2^{(r-1)h+i} \leq x < 2^{rh+i}.
\]
If $a_i \in A_i$ and $a_i \leq x$, then there is a set
\[
F \subseteq \{ i, h+i,\ldots, (r-1)h+i\}
\]
such that
\[
a_i = \sum_{f\in F} 2^f.
\]
The number of such sets $F$ is exactly $2^r$.  Since $0 \in A_i,$ we have
\[
A_i(x) \leq 2^r -1 < 2^r \leq 2^{1-i/h}x^{1/h} 
\]
and so
\bq
A(x) & = & A_0(x)+A_1(x)+\cdots + A_{h-1}(x) \\
& < & \left(\sum_{i=0}^{h-1} 2^{1-i/h}\right) x^{1/h} \\
& = & \left(\frac{1}{1-2^{-1/h}}\right)x^{1/h}.
\eq
Thus, $A$ is a thin basis of order $h$.  This completes the proof.
\end{proof}

For $h = 2,$ the Raikov-St\" ohr construction produces the thin basis
$A = A_0 \cup A_1$ of order 2, where
\[
A_0 = \{0,1,4,5,16,17,20,21,64,65,68,69,80,81,84,85,256,\ldots\}
\]
is the set of all finite sums of even powers of 2, and
\[
A_1 = \{0,2,8,10,32,34,40,42,128, 130, 136,138, 160,162,168,170,512,\ldots\}
\]
is the set of all finite sums of odd powers of 2.

\section{Construction of thin $g$-adic bases of order $h$}

\bl                      \label{thin:lemma:gw-adic}
Let $g \geq 2.$
Let $W$ be a nonempty set of nonnegative integers such that 
\[
W(x) = \theta x + O(1)
\]
for some $\theta \geq 0$ and all $x \geq 1$.
Let $\mathcal{F}(W)$ be the set of all finite subsets of $W$.
Let $A(W)$ be the set consisting of all integers of the form
\beq            \label{thin:gw-adic}
a = \sum_{w\in F} e_w g^w
\eeq
where $F\in \mathcal{F}(W)$ and $e_w \in \{0,1,\ldots,g-1\}$
for all $w\in F$.   Then
\[
x^{\theta}  \ll A(W)(x) \ll x^{\theta}
\]
for all sufficiently large $x$.
\el

\begin{proof}
The nonempty set $W$ is finite if and only if $\theta = 0$, and in this case $A(W)$ is also nonempty and finite, or, equivalently, $1 \ll A(W)(x) \ll 1$.  

Suppose that $\theta > 0$ and the set $W$ is infinite.
Let $W = \{w_i\}_{i=1}^{\infty}$, where $0 \leq w_1 < w_2 < w_3 < \cdots$.
Let $\delta = 0$ if $w_1 \geq 1$ and $\delta = 1$  if $w_1 = 0$.  
For $x\geq g^{w_1}$,  we choose the positive integer $k$ so that
\[
g^{w_k} \leq x < g^{w_{k+1}}.
\]
Then 
\[
w_k \leq \frac{\log x}{\log g} < w_{k+1}
\]
and
\[
k = W\left(\frac{\log x}{\log g} \right) + \delta 
= \frac{\theta \log x}{\log g}+ O(1) 
\]
where $W(x)$ is the counting function of the set $W$.

If $a \in A(W)$ and $a \leq x,$ then every power of $g$ that appears
with a nonzero coefficient in the $g$-adic 
representation~(\ref{thin:gw-adic}) of $a$ 
does not exceed $g^{w_k}$, and so $a$ can be written in the form
\[
a = \sum_{i=1}^k e_{w_i} g^{w_i},
\qquad\mbox{where $e_{w_i} \in \{0,1,\ldots,g-1\}$.}
\]
There are exactly $g^{k}$ integers of this form, and so
\[
A(W)(x) \leq g^{k}= g^{\frac{\theta \log x}{\log g}+O(1)} \ll x^{\theta}.
\]
Similarly, if $a$ is one of the $g^{k-1}-1$ positive integers
that can be represented in the form
\[
a = \sum_{i=0}^{k-1} e_{w_i} g^{w_i},
\]
then 
\[
a \leq  \sum_{i=0}^{k-1} (g-1) g^{w_i} \leq \sum_{j=0}^{w_{k-1}} (g-1) g^j
< g^{w_{k-1} +1}  \leq g^{w_k} \leq x
\]
and so
\[
A(W)(x) \geq g^{k-1}-1 \gg x^{\theta}.
\]
This completes the proof.
\end{proof}

\bt[Jia-Nathanson]         \label{thin:theorem:Jia-Nathanson}
Let $g \geq 2$ and $h \geq 2.$
Let $W_0, W_1,\ldots, W_{h-1}$ be nonempty sets of nonnegative integers such that
\[
\N_0 = W_0 \cup W_1 \cup \cdots \cup W_{h-1}
\]
and 
\[
W_i(x) = \theta_ix + O(1)
\]
where $0 \leq \theta_i \leq 1$ for $i = 0,1,\ldots, h-1$.
Let
\[
\theta = \max(\theta_0, \theta_1,\ldots,\theta_{h-1}).
\]
Let $A(W_0), A(W_1),\ldots, A(W_{h-1})$ be the sets of nonnegative integers constructed in Lemma~\ref{thin:lemma:gw-adic}.
The set
\[
A = A(W_0) \cup A(W_1) \cup \cdots \cup A(W_{h-1})
\]
is a basis of order $h$, and
\[
A(x) = O\left(x^{\theta}\right).
\]
In particular, if 
\[
W_i(x) = \frac{x}{h} + O(1)
\]
for $i = 0,1,\ldots, h-1$, then $A = A(W_0) \cup A(W_1) \cup \cdots \cup A(W_{h-1})$  is a thin basis of order $h$.
\et

Note that it is not necessary to assume that the sets $W_0,W_1,\ldots, W_{h-1}$ are pairwise disjoint.  

\begin{proof}
Every nonnegative integer $n$ has a $g$-adic representation 
of the form
\[
n = \sum_{w=0}^t e_wg^w,
\]
where $t \geq 0$ and $e_w \in \{0,1,\ldots,g-1\}$ for $w = 0,1,\ldots,t.$
We define the sets 
\bq
F_0 & = & \{ w\in \{0,1,\ldots,t\} : w\in W_0\} \\
F_1 & = & \{ w\in \{0,1,\ldots,t\} : w\in W_1\setminus W_0\} \\
F_2 & = & \{ w\in \{0,1,\ldots,t\} : w\in W_2\setminus (W_0 \cup W_1) \} \\
& \vdots & \\
F_{h-1} & = & \{ w\in \{0,1,\ldots,t\} : w\in W_{h-1}\setminus (W_0 \cup\cdots\cup W_{h-2}) \}.
\eq
Then $F_i \in \mathcal{F}(W_i)$ for all $i = 0,1,\ldots,h-1$.  Since $0 \in A(W_i)$ for $i = 0,1,\ldots, h-1$, we have 
\[
n = \sum_{w=0}^t e_wg^w = \sum_{i=0}^{h-1} \sum_{w\in F_i} e_w g^w
\in A(W_0)+\cdots + A(W_{h-1}) \in hA.
\]
Thus, $A$ is a basis of order $h$.

By Lemma~\ref{thin:lemma:gw-adic}, 
\[
A(W_i)(x) = O\left(x^{\theta_i}\right) = O\left(x^{\theta}\right)
\]
for all $i = 0,1,\ldots,h-1,$ and so
\[
A(W)(x) \leq \sum_{i=0}^{h-1} A(W_i)(x) = O\left(x^{\theta}\right).
\]
If $\theta_i = 1/h$ for all $i$, then $\theta = 1/h$ and $A$ is a thin basis.  This completes the proof.
\end{proof}

Consider the case when $W_i = \{w\in \N_0 : w \equiv i\pmod {h}\}$ for $i=0,1,\ldots,h-1$.  We shall compute an upper bound for the counting functions $A_i(x)$ and $A(x).$  For each $i$ and $x \geq g^i$, choose the positive integer $r$ such that 
\[
g^{(r-1)h+i} \leq x < g^{rh+i}.
\]
Then 
\[
A_i(x) \leq g^r -1 < g^r \leq g^{1-(i/h)}x^{1/h}
\]
and so
\[
A(x) = \sum_{i=0}^{h-1} A_i(x) < \sum_{i=0}^{h-1} g^{1-(i/h)}x^{1/h} = 
\frac{g-1}{1-g^{-1/h}} x^{1/h}. 
\]
Applying the mean value theorem to the function $f(x) = x^{1/h},$ we obtain $A(x) < ghx^{1/h}.$  In particular, if $g=2,$ we obtain 
$A(x)  < \frac{1}{1-2^{-1/h}} x^{1/h}< 2hx^{1/h}.$
This special case is the Raikov-St\" ohr construction.  
For $h=2$ the Raikov-St\" ohr basis $A = \{a_k\}_{k=1}^{\infty}$ with $a_k < a_{k+1}$ for $k \geq 1$  satisfies
\[
\frac{A(x)}{\sqrt{x}} < 2 + \sqrt{2} = 3.4142\ldots .
\]
Letting $x=a_k,$ we obtain 
\[
\frac{a_k}{k^2} > \frac{3-2\sqrt{2}}{2} = 0.0857 \ldots.
\]

If $A$ is a basis of order $h$, then the order of magnitude of the counting function
$A(x)$ must be at least $x^{1/h},$ and there exist thin bases,
such as the Raikov-St\" ohr bases and the Jia-Nathanson bases, with exactly this order of magnitude.
Two natural constants associated with thin bases of order $h$ are
\[
\alpha_h = \inf_{ \substack{ A\subseteq \N_0 \\ hA=\N_0}}
\liminf_{x\rightarrow\infty}\frac{A(x)}{x^{1/h}}
\]
and
\[
\beta_h = \inf_{ \substack{ A\subseteq \N_0 \\ hA=\N_0} }
\limsup_{x\rightarrow\infty}\frac{A(x)}{x^{1/h}}
\]
St\" ohr~\cite{stoh55} proved the following lower bound for $\beta_h$.

\bt[St\" ohr]   \label{thin:theorem:Stohr}
\[
\beta_h \geq \frac{\sqrt[h]{h!}}{\Gamma(1+1/h)}
\]
where $\Gamma(x)$ is the Gamma function.    
\et

In particular,  $\limsup_{x\rightarrow\infty} A(x)/\sqrt{x} \geq \sqrt{8/\pi}$ for every basis $A$ of order 2.  
   
\boprob
Compute the numbers $\alpha_h$ and $\beta_h$
for all $h \geq 2.$   
\eoprob

This is an old unsolved problem in additive number theory.  Even the numbers $\alpha_2$ and $\beta_2$ are unknown.

\section{Asymptotically polynomial bases}
Let $h \geq 2$, and let  $A = \{a_k\}_{k=1}^{\infty}$ be a set of nonnegative integers with $a_1 = 0$ and $a_k < a_{k+1}$ for all $k \geq 1$.  If $A$ is a basis of order $h$, then there is a real number $\lambda_2$ such that $ a_k \leq \lambda_2 k^h$ for all $k$ (Theorem~\ref{thin:theorem:thin1}).  The basis $A$ is called \emph{thin} if there is also a number $\lambda_1 > 0$ such that $a_k \geq \lambda_1k^h$ for all $k$.  Thus, if $A$ is a thin basis of order $h$, then there exist positive real numbers $\lambda_1$ and $\lambda_2$ such that 
\[
\lambda_1 \leq \frac{a_k}{k^h} \leq \lambda_2
\]
for all $k$.  In Theorems~\ref{thin:theorem:Raikov-Stohr} and \ref{thin:theorem:Jia-Nathanson} we constructed examples of thin bases of order $h$ for all $h \geq 2$.  

The sequence  $A = \{a_k\}_{k=0}^{\infty}$ is called \index{asymptotically polynomial} \emph{asymptotically polynomial of degree $d$} if there is a real number $\lambda > 0$ such that $a_k \sim \lambda k^d$ as $k \rightarrow \infty.$  If $A$ is a basis of order $h$ and if $A$ is also asymptotically polynomial of degree $d$, then $d \leq h$.  
We shall describe a beautiful construction of J. W. S. Cassels of a family of additive bases of order $h$ that are asymptotically polynomial of degree $h$.   The key to the construction is the following result, which allows us to embed a sequence of nonnegative integers with regular growth into a sequence of nonnegative integers  with asymptotically polynomial growth.

\bt         \label{Cassels:theorem:casimbed}
Let $h\geq 2$ and let $A = \{a_k\}_{k=1}^{\infty}$ 
be a sequence of nonnegative integers such that 
\[
\liminf_{k\rightarrow\infty} \frac{a_{k+1}-a_k}{a_k^{(h-1)/h}} = \alpha > 0
\]
For every real number $\gamma$ with $0 < \gamma < \alpha$,  there exists a sequence $C = \{c_k\}_{k=0}^{\infty}$ of nonnegative integers such that $C$ is a supersequence of $A$ and 
\[
c_k = \left( \frac{\gamma k}{h}\right)^h + O\left(k^{h-1}\right).
\]
\et

\begin{proof}
Let $B = \{b_k\}_{k=1}^{\infty}$ be a strictly increasing sequence 
of nonnegative  integers such that 
\[
b_k = \left(\frac{\gamma k}{h}\right)^h  + O(k^{h-2}).
\]
Since $h \geq 2$ and $b_k =  \left(\gamma k/h\right)^h \left(1 + O\left(k^{-2}\right)\right)$, we have 
\begin{align*}
\frac{b_{k+1}-b_k}{b_k^{(h-1)/h}} 
& = \frac{  \left(\frac{\gamma}{h}\right)^h \left( (k+1)^h-k^h+O(k^{h-2})\right)}{\left(\frac{\gamma k}{h}\right)^{h-1}  \left(1+O\left(k^{-2}\right)  \right)^{(h-1)/h}} \\
& = \frac{ \gamma \left( hk^{h-1}+ O\left(k^{h-2}\right)\right)}{hk^{h-1}  \left(1+O\left(k^{-2}\right)  \right)^{(h-1)/h}} \\
& =  \frac{\gamma \left( 1 +O\left(k^{-1} \right) \right)}
{\left(1+O\left(k^{-2}\right)  \right)^{(h-1)/h}} \\
& =  \gamma(1+ o(1))
\end{align*}
and so
\[
\lim_{k\rightarrow \infty} \frac{b_{k+1}-b_k}{b_k^{(h-1)/h}} = \gamma.
\]

Suppose there exist infinitely many $k$ such that, for some integer $m = m(k)$,
\[
b_k < a_m < a_{m+1} \leq b_{k+1}.
\]
The inequality
\[
\frac{b_{k+1} - b_k}{b_k^{(h-1)/h}} > \frac{a_{m+1} - a_m}{b_k^{(h-1)/h}} 
 > \frac{a_{m+1} - a_m}{a_m^{(h-1)/h}}
\]
implies that
\[
\gamma = \lim_{k\rightarrow\infty}\frac{b_{k+1} - b_k}{b_k^{(h-1)/h}} 
\geq \liminf_{m\rightarrow\infty}\ \frac{a_{m+1} - a_m}{a_m^{(h-1)/h}} \geq \alpha > \gamma
\]
which is impossible.  Therefore, there exists an integer $K$ such that, for every integer $k \geq K,$
the interval $(b_k,b_{k+1}]$ contains at most one element of $A$.

Choose the integer $L$ such that
\[
a_L \leq b_K < a_{L+1}.
\]
We define the sequence $C =\{c_k\}_{k=0}^{\infty}$ as follows:
Let $c_k = a_k$ for $k = 1,2,\ldots, L.$
For $i \geq 1,$ we choose $c_{L+i} \in (b_{K+i-1}, b_{K+i}]$ as follows:  If the interval $(b_{K+i-1}, b_{K+i}]$ contains the element $a_{\ell}$ from the sequence $A$, then $c_{L+i} = a_{\ell}$.  Otherwise, let $c_{L+i} = b_{K+i}$.   Since the interval $(b_{K+i-1}, b_{K+i}]$
contains at most one element of $A$ for all $i \geq 1$, and since every element $a_k$ of $A$ with $k > L$ 
is contained in some interval of the form $(b_{K+i-1},b_{K+i}]$ with $i \geq 1$, it follows that $A$ is a subsequence of $C$.
Moreover, for every $k \geq L+1,$ 
\[
b_{k-L+K-1} < c_{k} \leq b_{k-L+K}.
\]
Since
\[
b_{k-L+K} = \left(\frac{\gamma}{h}\right)^h (k-L+K)^h + O(k^{h-2})
= \left( \frac{\gamma k}{h}\right)^h  + O\left(k^{h-1}\right)
\]
and, similarly, $b_{k-L+K-1} = \left( \gamma k/h\right)^h  + O\left(k^{h-1}\right)$, it follows that
\[
c_k = \left(\frac{\gamma k}{h}\right)^h  + O\left(k^{h-1}\right).
\]
This completes the proof.
\end{proof}

\section{Bases of order 2}

In this section we describe Cassels' construction in the case $h = 2$.
We need the following convergence result.

\bl              \label{Cassels:lemma:converge}
Let $0 < \alpha < 1$.  If $\{q_k\}_{k=1}^{\infty}$ is a sequence of positive integers such that
\[
\lim_{k\rightarrow\infty} \frac{q_{k-1}}{q_{k}} = \alpha
\]
then
\[
\lim_{k\rightarrow\infty} \frac{q_1+q_2 + \cdots + q_k}{q_k} = \frac{1}{1-\alpha}.
\]
\el

\begin{proof}
For every nonnegative integer $j$ we have                
\beq                \label{Cassels:ratio}
\lim_{k\rightarrow\infty} \frac{q_{k-j}}{q_{k}} 
= \lim_{k\rightarrow\infty} \prod_{i=0}^{j-1} \frac{q_{k-i-1}}{q_{k-i}}
= \alpha^j.
\eeq
Let $\beta$ be a real number such that $\alpha < \beta < 1.$
For every $\varepsilon > 0$ there exists a number $K = K(\beta,\varepsilon)$
such that 
\beq                 \label{Cassels:qk}
\frac{q_{k-1}}{q_{k}} < \beta \qquad\mbox{for all $k \geq K$}
\eeq
and
\beq           \label{Cassels:beta}
\beta^K < \frac{(1-\beta)\varepsilon }{4}.
\eeq
If $k \geq K$ and $k-K = r$, then 
\[
q_k  > \beta^{-1}q_{k-1} > \beta^{-2}q_{k-2} > \cdots > \beta^{-r}q_{k-r}
= \beta^{K-k}q_{K} = c\beta^{-k}
\]
where $c = \beta^{K}q_{K} > 0$, and so 
\beq           \label{Cassels:inf}
\lim_{k\rightarrow\infty}q_k = \infty.
\eeq
If $0 \leq j \leq k-K+1$, then inequality~\eqref{Cassels:qk} implies
\[
\frac{q_{k-j}}{q_{k}} = \prod_{i=0}^{j-1} \frac{q_{k-i-1}}{q_{k-i}} < \beta^j.
\]
For $k \geq 2K$ we obtain
\bq
\lefteqn{\left|\frac{q_1+q_2 + \cdots + q_k}{q_k} - \frac{1}{1-\alpha}\right|
 =  \left|\sum_{j=0}^{k-1} \frac{q_{k-j}}{q_k} - \sum_{j=0}^{\infty} \alpha^j \right|} \\
& \leq & \sum_{j=0}^{K-1}\left| \frac{q_{k-j}}{q_k} -  \alpha^j \right| 
+ \sum_{j=K}^{k-K+1} \frac{q_{k-j}}{q_k} 
+ \sum_{j=k-K+2}^{k-1} \frac{q_{k-j}}{q_k} 
+ \sum_{j=K}^{\infty} \alpha^j \\
& < & \sum_{j=0}^{K-1}\left| \frac{q_{k-j}}{q_k} -  \alpha^j \right| 
+ \sum_{j=K}^{k-K+1} \beta^j
+ \sum_{j=1}^{K-2} \frac{q_{j}}{q_k} 
+ \sum_{j=K}^{\infty} \beta^j \\
& < & \sum_{j=0}^{K-1}\left| \frac{q_{k-j}}{q_k} -  \alpha^j \right| 
+ \sum_{j=1}^{K-2} \frac{q_{j}}{q_k} 
+ \frac{2\beta^K}{1-\beta}
\eq
It follows from~(\ref{Cassels:ratio}),~(\ref{Cassels:inf}), 
and~\eqref{Cassels:beta}
that for  $j = 0,1,\ldots,K-1$ and all sufficiently large $k$ 
\[
\left| \frac{q_{k-j}}{q_k} -  \alpha^j \right| < \frac{\varepsilon}{4K}
\]
and
\[
\frac{q_j}{q_k} < \frac{\varepsilon}{4K}
\]
and so
\[
\left|\frac{q_1+q_2 + \cdots + q_k}{q_k} - \frac{1}{1-\alpha}\right| < \varepsilon.
\]
This completes the proof.
\end{proof}

\bt              \label{Cassels:theorem:APbasis} 
Let $\{q_i\}_{i=1}^{\infty}$ and $\{m_i\}_{i=1}^{\infty}$
be sequences of positive integers such that
\beq  \label{Cassels:AP1}
q_1 = 1
\eeq
and, for all $i \geq 2,$
\beq  \label{Cassels:AP2}
(q_{i-1},q_i) = (q_{i-1},q_{i+1}) = 1
\eeq
\beq  \label{Cassels:AP3}
m_{i-1} \geq q_{i}+q_{i+1}-2
\eeq
and
\beq  \label{Cassels:AP4}
m_{i+1}q_{i+1} \geq m_{i}q_{i} + m_{i-1}q_{i-1}.
\eeq
Define the sequences $\{Q_k\}_{k=1}^{\infty}$ of nonnegative integers
and $\{A_k\}_{k=1}^{\infty}$ of finite arithmetic progressions 
of nonnegative integers by
\[
Q_k = \sum_{i=1}^{k-1} m_iq_i 
\]
and
\[
A_k = Q_k + q_k\ast [0,m_k].
\]
Let
\[
A = \bigcup_{k=1}^{\infty} A_k  = \{a_n\}_{n=0}^{\infty} 
\]
where $a_0 = 0 < a_1 < a_2 < \cdots$.  Then $A$ is a basis of order 2, and, for every positive integer $K$, the set
$\bigcup_{k=K}^{\infty} A_k$
is an asymptotic basis of order 2.

Let  $A(x)$ be the counting function of the set $A$, and let 
$M_k = \sum_{i=1}^{k-1} m_i$ for $k \geq 1$.  
If $M_k \leq n \leq M_{k+1},$ then
\beq   \label{Cassels:cas1}
a_n = Q_k + (n-M_k)q_k.
\eeq
If $Q_k \leq x \leq Q_{k+1},$ then
\beq   \label{Cassels:cas2}
A(x) = M_k + \left[\frac{x-Q_k}{q_k}\right].
\eeq
\et

\begin{proof}
Since $Q_{k+1} - Q_k = m_kq_k,$ it follows that 
\[
\{Q_k, Q_{k+1}\} \subseteq A_k \subseteq [Q_k,Q_{k+1}]
\]
and
\[
A_k = Q_{k+1} - q_k\ast [0,m_k].
\]
Also, $Q_1 = 0$, $Q_2 = m_1q_1 = m_1$,  
and $A_1 = [0,m_1],$ hence
\[
 [2Q_1,2Q_2] = [0,2m_1] =2A_1.
\]
We shall prove that
\beq  \label{Cassels:BasisEstimate}
[2Q_k, 2Q_{k+1}] \subseteq 
A_{k-1} + \left(A_k \cup A_{k+1}\right)
\subseteq 2\left( A_{k-1} \cup A_k \cup A_{k+1} \right)
\eeq
for all $k \geq 2.$

Let $n \in [2Q_k, 2Q_{k+1}].$  
There are two cases.
In the first case we have 
\beq         \label{Cassels:case1}
2Q_k \leq n \leq Q_k + Q_{k+1} - (q_k-1)q_{k-1}.
\eeq
Since $(q_k,q_{k-1}) = 1$, there is a unique integer $r$ such that
\[
n \equiv 2Q_k - rq_{k-1} \pmod{q_k}
\]
and, by~\eqref{Cassels:AP3},
\beq         \label{Cassels:xqk1}
0 \leq r \leq q_k -1 \leq m_{k-1}.
\eeq
Then $Q_k - rq_{k-1} \in A_{k-1}.$
There is a unique integer $s$ such that
\[
sq_k = n - 2Q_k + rq_{k-1}.
\]
It follows from~(\ref{Cassels:case1}) and~(\ref{Cassels:xqk1}) that
\[
0 \leq n - 2Q_k + rq_{k-1} \leq  Q_{k+1} - Q_k = m_kq_k,
\]
and so
\[
0 \leq s \leq m_k. 
\]
Therefore, $Q_k + sq_k \in A_k$ and 
\[
n  = \left(Q_k-rq_{k-1}\right) + \left(Q_k + sq_k\right) \in A_{k-1} + A_k.
\]

In the second case we have 
\beq         \label{Cassels:case2}
Q_k + Q_{k+1} - (q_k-1)q_{k-1} +1 \leq n \leq 2Q_{k+1}.
\eeq
The set $R = [q_k -1, q_k + q_{k+1}-2]$ is a complete set of representatives of the congruence classes modulo $q_{k+1}$.
Since $(q_{k-1},q_{k+1}) = 1$,  it follows that there is a unique integer $r \in R$ such that
\[
n \equiv Q_k + Q _{k+1} - rq_{k-1} \pmod{q_{k+1}}.
\]
Inequality~\eqref{Cassels:AP3} implies that 
\beq        \label{Cassels:xqk2}
0 \leq q_k -1 \leq r \leq q_k + q_{k+1}-2 \leq m_{k-1}
\eeq
and so  $Q_k - rq_{k-1} \in A_{k-1}$.
There is a unique integer $t$ such that
\[
tq_{k+1} = n - Q_k - Q _{k+1} + rq_{k-1},
\]
Inequalities~\eqref{Cassels:case2},~\eqref{Cassels:xqk2}, and~\eqref{Cassels:AP4} imply that
\[
tq_{k+1} \geq  (r - q_k + 1)q_{k-1} +1 \geq 1
\]
and
\[
tq_{k+1} \leq Q_{k+1} - Q_k + rq_{k-1} \leq m_kq_k + m_{k-1}q_{k-1} \leq m_{k+1}q_{k+1},
\]
and so
\[
1 \leq t \leq m_{k+1}.
\]
Therefore, $Q_{k+1} + tq_{k+1}\in A_{k+1}$ and 
\[
n = (Q_k - rq_{k-1}) + (Q_{k+1} + tq_{k+1}) \in A_{k-1}+A_{k+1}.
\]
This proves~\eqref{Cassels:BasisEstimate}.  
It follows that $\bigcup_{k=1}^{\infty} A_k$ is a basis of order 2.  
Moreover, for every positive integer $K$,
\[
\left[2Q_{K+1},\infty\right) \subseteq 2 \left(\bigcup_{k=K}^{\infty} A_k\right)
\]
and so $\bigcup_{k=K}^{\infty} A_k$ is an asymptotic basis of order 2.

Let $A = \{a_n\}_{n=0}^{\infty},$ where $a_0 = 0 < a_1 < a_2 < \cdots$,
and let $A(x)$ be the counting function of the set $A$.  
Formulas~(\ref{Cassels:cas1}) and~(\ref{Cassels:cas2}) are immediate consequences of the construction of the set $A$. 
This completes the proof.
\end{proof}

\bt              \label{Cassels:theorem:C2} 
Let $0 < \alpha < 1$ and let $\{q_i\}_{i=1}^{\infty}$ be a sequence of positive integers 
with $q_1 = 1$ such that, for all $i \geq 2,$
\beq   \label{Cassels:theorem:C2-a} 
(q_{i-1},q_i) = (q_{i-1},q_{i+1}) = 1
\eeq
\beq   \label{Cassels:theorem:C2-b} 
q_{i+1} (q_{i+2} + q_{i+3}) \geq q_{i}(q_{i+1}+q_{i+2}) + q_{i-1}(q_{i}+q_{i+1})
\eeq
and
\beq   \label{Cassels:theorem:C2-c} 
\lim_{i\rightarrow\infty} \frac{q_{i-1}}{q_i} = \alpha.
\eeq
Define the sequences $\{Q_k\}_{k=1}^{\infty}$ of nonnegative integers
and $\{A_k\}_{k=1}^{\infty}$ of finite arithmetic progressions 
of nonnegative integers by
\[
Q_k = \sum_{i=1}^{k-1} q_i(q_{i+1}+q_{i+2})  
\]
and
\[
A_k = Q_k + q_k\ast [0,q_{k+1}+q_{k+2} ].
\]
Let 
\[
A = \bigcup_{k=1}^{\infty} A_k = \{a_n\}_{n=0}^{\infty},
\]
where $a_0 = 0 < a_1 < a_2 < \cdots$.
Then $A$ is a basis of order 2 such that 
\[
\liminf_{k\rightarrow\infty} \frac{a_{n+1}-a_n}{n} 
\geq \frac{\alpha^2(1-\alpha)}{1+\alpha} > 0.
\]
\et

Note that the sequence $\{q_i\}_{i=1}^{\infty}$ of Fibonacci numbers defined by $q_1 = q_2 = 1$ and $q_{i+2} = q_{i+1}+q_i$ for $i \geq 1$ satisfies the conditions of Theorem~\ref{Cassels:theorem:C2} with $\alpha = (\sqrt{5}-1)/2$.

\begin{proof}
For every  integer $i \geq 1$ we define the positive integer $m_i = q_{i+1}+q_{i+2}$.  Inequality~\eqref{Cassels:theorem:C2-b} implies that the sequence $\{m_i\}_{i=1}^{\infty}$ satisfies the hypotheses of Theorem~\ref{Cassels:theorem:APbasis}, and so $A$ is a basis of order 2.
For $k \geq 1$ we define 
\[
M_k = \sum_{i=1}^{k-1} m_i= \sum_{i=1}^{k-1} (q_{i+1}+q_{i+2}).
\]
Then $\{M_k\}_{k=1}^{\infty}$ is a strictly increasing sequence of positive integers.  
For every positive integer $n$ there is a unique integer $k$ such that
\[
M_k \leq n < M_{k+1}.
\]
By ~\eqref{Cassels:cas1} we have 
\[
a_n = Q_k + (n-M_k)q_k
\]
and so
\[
a_{n+1}-a_n = q_k,
\]
hence
\[
\frac{a_{n+1}-a_n}{n} = \frac{q_k}{n} > \frac{q_k}{M_{k+1}}. 
\]
Condition~\eqref{Cassels:theorem:C2-c} implies that $\lim_{k\rightarrow\infty} q_k = \infty$.  Since
\bq
\frac{M_{k+1}}{q_k} 
& = & \sum_{i=1}^{k} \frac{q_{i+1}+q_{i+2}}{q_k} \\
& = & \sum_{i=2}^{k+1} \frac{q_{i}}{q_k} + \sum_{i=3}^{k+2} \frac{q_{i}}{q_k} \\
& = & 2\sum_{i=1}^{k} \frac{q_{i}}{q_k} + 2\frac{q_{k+1}}{q_k}
+ \frac{q_{k+2}}{q_k} - 2\frac{q_{1}}{q_k} - \frac{q_{2}}{q_k},
\eq
it follows from Lemma~\ref{Cassels:lemma:converge} that
\[
\lim_{k\rightarrow\infty} \frac{M_{k+1}}{q_k} 
= \frac{2}{1-\alpha} + \frac{2}{\alpha} + \frac{1}{\alpha^2}
= \frac{1+\alpha}{\alpha^2(1-\alpha)}.
\]
Therefore,
\[
\liminf_{k\rightarrow\infty} \frac{a_{n+1}-a_n}{n} 
\geq \lim_{k\rightarrow\infty} \frac{q_k} {M_{k+1}}
= \frac{\alpha^2(1-\alpha)}{1+\alpha} > 0.
\]
This completes the proof.
\end{proof}

\bt[Cassels]
There  exist a basis $C = \{c_n\}_{n=0}^{\infty}$ of order 2 and a real number $\lambda > 0$ such that  $c_n = \lambda n^2 + O(n).$
\et

\begin{proof}
By Theorem~\ref{Cassels:theorem:C2}, there exists a basis $A = \{a_n\}_{n=0}^{\infty}$ of order 2 such that $\liminf_{n\rightarrow\infty} (a_{n+1}-a_n)/n  > 0$.   Applying Theorem~\ref {thin:theorem:thin1} with $h=2$, 
we see that $a_n \ll n^2$ and so $\liminf_{n\rightarrow\infty} (a_{n+1}-a_n)/a_n^{1/2}  > 0$.  Applying Theorem~\ref{Cassels:theorem:casimbed} with $h=2$, we obtain a sequence $C = \{c_n\}_{n=0}^{\infty}$ of nonnegative integers and a positive real number $\lambda$ such that $C$ is a supersequence of $A$ and $c_n = \gamma n^2 + O\left(n \right).$  This completes the proof.
\end{proof}

\section{Bases of order $h \geq 3$}

We start with Cassels' construction of a finite set $C$\ of integers such that the elements of $C$\ are widely spaced and $C$\ is a basis of order $h$ for a long interval of integers.  The construction uses a perturbation  of the $g$-adic representation.  

\bl       \label{Cassels:lemma:IntervalBasis}
Let $h \geq 3.$  Let  $v$ and $L$ be positive integers with $L \geq h.$  Define
\[
g = 2^{h+1}v.
\]
Let $C = C(v,L)$ denote the finite set consisting of the following integers:
\[
\ba{ll}
g^h + eg^{h-1} + 2vg^{h-2} + e & \mbox{for $0 \leq e < g$,} \\
(i+1)g^h + eg^{h-1} + eg^i & \mbox{for $0 \leq i \leq h-3$ and $0 \leq e < g$,} \\
(h-1)g^h + (4vq+r)g^{h-1} 
+ (4vq+r)g^{h-2} & \mbox{for $0 \leq q < 2^{h-1}$ and $0 \leq r < 2v$,} \\
hg^h + \ell g^{h-1}  & \mbox{for $0 \leq \ell < Lg$.}
\ea
\]
Then
\benum
\item[(i)]
The $h$-fold sumset $hC$ contains every integer $n$ in the interval
\[
\left[ \left(\frac{h^2 + 3h -2}{2}  \right)g^h , \left(\frac{h(h+1)}{2} + L\right)g^h\right).
\]
\item[(ii)]
If $c\in C$, then 
\[
g^h \leq c < (h+L)g^h.
\]
If $c \geq hg^h$, then $c \equiv 0 \pmod{g^{h-1}}$.
\item[(iii)]
If $c,c' \in C$ and $c \neq c'$, then 
\[
|c - c'| \geq vg^{h-2} - g.
\]
\item[(iv)]
If $c \in C$ and $y$ is any integer such that
\[
y \equiv -vg^{h-2} \pmod{4vg^{h-2}}
\]
then
\[
|c-y| \geq vg^{h-2}-g.
\]
\eenum
\el

\begin{proof}
(i)  Every nonnegative integer $n$ has a unique $g$-adic representation
in the form
\beq  \label{Cassels:g-adicRep}
n =  e_{h-1}g^{h-1} + e_{h-2}g^{h-2} +\cdots + e_1g + e_0 
\eeq
where $e_{h-1} \geq 0$ and 
\[
0 \leq e_j < g \qquad\text{ for $j = 0,1,\ldots, h-2$}.
\]
If $n$ satisfies the inequality
\[
\left( \frac{h^2 + 3h -2}{2}\right)g^h \leq n < \left(\frac{h(h+1)}{2} + L\right)g^h
\]
then $e_{h-1}$ satisfies the inequality 
\beq    \label{Cassels:ehineq}
\left(  \frac{h^2 + 3h -2}{2} \right) g \leq e_{h-1} <  \left(\frac{h(h+1)}{2} + L\right) g.
\eeq
The digit $e_{h-2}$ satisfies the inequality $0 \leq e_{h-2} < g = 4v2^{h-1}.$
There are two cases, which depend on the remainder 
of $e_{h-2}$ when divided by $4v$.  

In the first case, we have
\[
e_{h-2} = 4vq + r \qquad\mbox{with $0 \leq q < 2^{h-1}$ and $0 \leq r < 2v$}
\]
Rearranging the $g$-adic representation~\eqref{Cassels:g-adicRep}, 
we obtain
\begin{multline}    \label{Cassels:Repn1}
n  =   \left( (h-1)g^h + (4vq + r )g^{h-1} + (4vq + r )g^{h-2} \right)  + \\
+ \sum_{i=0}^{h-3} \left( (i+1)g^h + e_ig^{h-1} + e_ig^i \right) + \left( hg^h + \ell g^{h-1} \right)
\end{multline}
where 
\[
\ell = e_{h-1} - \sum_{i=0}^{h-2}e_i -\frac{h(h+1)g}{2}.
\]
Inequality~\eqref{Cassels:ehineq} implies that
\[
\ell \geq \left( \frac{h^2 + 3h -2}{2}\right) g - (h-1)(g-1) -\frac{h(h+1)g}{2} = h-1 > 0
\]
and 
\[
\ell < \left(\frac{h(h+1)}{2} + L\right) g - \frac{h(h+1)g}{2} = Lg 
\]
and so $hg^h + \ell g^{h-1} \in C$.   Thus,~\eqref{Cassels:Repn1} is a representation of $n$ as the sum of $h$ elements of $C$, that is, $n \in hC$.  

In the second case, we have
\[
e_{h-2} = 4vq + r + 2v \qquad\mbox{with $0 \leq q < 2^{h-1}$ and $0 \leq r < 2v$}.
\]
From the $g$-adic representation~\eqref{Cassels:g-adicRep}, we obtain 
\begin{multline}   \label{Cassels:Repn2}
n   =  \left((h-1)g^h + (4vq+r)g^{h-1} + (4vq+r)g^{h-2}\right) + \\
       + \sum_{i=1}^{h-3} \left( (i+1)g^h + e_ig^{h-1} + e_ig^i \right)  + \\
        + \left(g^h + e_0g^{h-1} + 2vg^{h-2} + e_0 \right) 
        +  \left( hg^h + \ell g^{h-1} \right)  
\end{multline}
where 
\[
\ell = e_{h-1} - (e_{h-2} - 2v) - \sum_{i=0}^{h-3} e_i  - \left(\frac{h(h+1)}{2}\right)g.
\]
As in the first case, inequality~\eqref{Cassels:ehineq} implies that
$0 < h-1 \leq \ell < Lg$ and so $hg^h + \ell g^{h-1} \in C$.   Thus,~\eqref{Cassels:Repn2} is a representation of $n$ as the sum of $h$ elements of $C$, that is, $n \in hC$.  
This proves~(i).

To prove~(ii), we observe that the smallest element of $C$\ is $g^h$ and the largest is $hg^h+ (Lg-1)g^{h-1} < (h + L)g^h$.  If $c\in C$ and  $c \geq hg^h$, then $c = hg^h + \ell g^{h-1}$ for some nonnegative integer  $\ell < Lg$, hence $c \equiv 0 \pmod{g^{h-1}}$.  

To prove~(iii), we assert that every integer $c \in C$ satisfies an inequality of the form 
\beq  \label{Cassels:c-bound}
4svg^{h-2} \leq c < (4s+2)vg^{h-2} + g
\eeq
for some nonnegative integer $s$.  There are four cases to check.

If $c = g^h+eg^{h-1} + 2vg^{h-2}+e$ with $0 \leq e < g$, then we choose 
$s = 2^{h-1}(g+e)$.  Since 
\[
4svg^{h-2} = g^h+eg^{h-1}
\]
and
\[
(4s+2)vg^{h-2} +g = g^h+eg^{h-1} + 2vg^{h-2} + g
\]
it follows that $c$ satisfies~\eqref{Cassels:c-bound}.

If $c = (i+1)g^h+eg^{h-1} + eg^i$ with $0 \leq e < g$ and $0 \leq i \leq h-3$, then $c$ satisfies~\eqref{Cassels:c-bound} with $s = 2^{h-1}((i+1)g+e)$.

If $c = (h-1)g^h+(4vq+r)g^{h-1} + (4vq+r)g^{h-2}$ with $0 \leq q < 2^{h-1}$ and $0 \leq r < 2v$, then $c$ satisfies~\eqref{Cassels:c-bound} with $s = 2^{h-1}((h-1)g+ 4vq+r) + q$.

If $c = hg^h+ \ell g^{h-1}$ with $0 \leq \ell < Lg$, then $c$ satisfies~\eqref{Cassels:c-bound} with $s = 2^{h-1}(hg+\ell)$.

This proves ~\eqref{Cassels:c-bound}.  It follows that the distance between elements of $C$ that satisfy inequality~\eqref{Cassels:c-bound} for different values of $s$ is at least $2vg^{h-2} - g$.  
If $c$ and $c'$ are distinct elements of $C$ that satisfy inequality~\eqref{Cassels:c-bound} for the same value of $s$, and if $c' < c$, then we must have 
\[
0 < c-c' < 2vg^{h-2} + g.
\]
This can happen only if $c = g^h+eg^{h-1} + 2vg^{h-2}+e$ and $c' = g^h+eg^{h-1} + e$ with $0 \leq e < g$, and so  $c-c' = 2vg^{h-2}$.  This proves~(iii).

Finally, to prove~(iv), we observe that if $y \equiv -vg^{h-2}\pmod{4vg^{h-2}}$, then $y = 4s'vg^{h-2} - vg^{h-2}$ for some integer $s'$, and the distance between $y$ and any integer satisfying an inequality of the form~\eqref{Cassels:c-bound} is at least $vg^{h-2} - g$.  This completes the proof of the Lemma.
\end{proof}

\bl  \label{Cassels:lemma:pjuj}
For $h \geq 3$, let $v_i = 2^i$ and $g_i = 2^{h+1}v_i = 2^{i+h+1}$ for $i = 0,1,2,\ldots.$  Then
\[
p_j  = \sum_{i=0}^j v_ig_i^{h-2} < g_j^h.
\]
\el

\begin{proof}
We compute $p_j$ explicitly as follows:
\begin{align*}
p_j & = \sum_{i=0}^j v_ig_i^{h-2}  = \sum_{i=0}^j 2^i \left(2^{i+h+1}\right)^{h-2}  =  2^{(h-2)(h+1)}\sum_{i=0}^j 2^{(h-1)i}\\
& = 2^{(h-2)(h+1)}\left( \frac{2^{(h-1)(j+1)}-1}{2^{h-1}-1} \right)  = \frac{2^{h^2+hj-j-3}-2^{h^2-h-2}}{2^{h-1}-1} \\
& < 2^{h(j + h +1)} = g_j^h
\end{align*}
because, for $h \geq 3$, 
\begin{align*}
2^{h^2 + hj - j - 3} + 2^{h^2+hj+h} & < 2^{h^2+hj+h +1}  < 2^{h^2+hj+2h-1} \\
&  < 2^{h^2+hj+2h-1} + 2^{h^2-h-2}.
\end{align*}
\end{proof}

\bt      \label{Cassels:theorem:main3a}
Let $h \geq 3$.   There exists a strictly increasing sequence $A = \{a_k\}_{k=1}^{\infty}$ of nonnegative integers such that $A$ is a basis of order $h$ and 
\[
\liminf_{k\rightarrow \infty}  \frac{a_{k+1}-a_k}{a_k^{(h-1)/h}} \geq \frac{1}{2^{3h-1}}.
\]
\et

\begin{proof}
Let 
\[
A(-1) = \left[ 0,2^{h^2+2h} \right].
\]
We define
\[
L = 2^{2h}-h-1
\]
and, for $i = 0,1,2,\ldots$,
\[
v_i = 2^i 
\]
\[
 g_i = 2^{h+1}v_i = 2^{i+h+1}
\]
and
\[
p_j = \sum_{i=0}^j v_ig_i^{h-2}.
\]
For $j=0,1,2,\ldots$, let
\[
A(j) = p_j + C(v_j,L)
\]
where $C(v_j,L)$ is the finite set of positive integers constructed in Lemma~\ref{Cassels:lemma:IntervalBasis}.  We begin by proving that
\[
A = \bigcup_{j=-1}^{\infty} A(j)
\]
is a basis of order $h$.  

First, we observe that 
\[
I(-1) =  \left[ 0, h2^{h^2+2h} \right] = hA(-1) \subseteq hA
\]
and, by Lemma~\ref{Cassels:lemma:IntervalBasis},
\[
I(j) = \left[ hp_j +   \left(\frac{h^2 + 3h -2)}{2} \right)g_j^h ,  hp_j +  \left(\frac{h(h+1)}{2} + L\right)g_j^h \right) \subseteq hA(j) 
\]
for $j = 0,1,2,\ldots$.  
Since $h^2+3h-2 \leq 2^{h+1}$ for $h \geq 3$, it follows that 
\begin{align*}
hp_0 +   \left(\frac{h^2 + 3h - 2)}{2}  \right)g_0^h 
& = h2^{(h+1)(h-2)} + \left(\frac{h^2+3h-2}{2}  \right)2^{h(h+1)}\\
& \leq h2^{h^2-h-2} + 2^{h^2+2h}\\
& \leq h2^{h^2+2h}
\end{align*} 
and so the intervals $I(-1)$ and $I(0)$ overlap.
Similarly, for $j \geq 0$ the intervals $I(j)$ and $I(j+1)$ overlap if 
\beq  \label{Cassels:IntervalIneq2}
hp_{j+1} +   \left(\frac{h^2+3h-2}{2}  \right)g_{j+1}^h 
\leq hp_j + \left( \frac{h(h+1)}{2} + L\right)g_j^h.
\eeq
Since $v_{j+1}=2v_j$ and $g_{j+1}=2g_j$, we have
\[
p_{j+1}- p_j = v_{j+1}g_{j+1}^{h-2} = 2^{h+j-1}g_j^{h-2}
= \frac{g_j^h}{2^{h+j+3}}.
\] 
Rearranging inequality~\eqref{Cassels:IntervalIneq2} and dividing by $g_j^h$, we see that it suffices to prove that    
\[
 \frac{h}{2^{h+j+3}} + \left( \frac{h^2+3h-2}{2}\right)2^h \leq  \frac{(h-2)(h+1)}{2} +  2^{2h}. 
\]
This follows immediately from the inequalities
$h^2+3h-2 \leq 2^{h+1}$ and
\[
\frac{h}{2^{h+j+3}} \leq 2 \leq \frac{(h-2)(h+1)}{2}
\]
for $j \geq 0$ and $h \geq 3$.  Thus, the set $A$ is a basis of order $h$.  

Next, we show that the elements of $A$ are widely spaced.  
Let $a,a' \in A$ with $a'\neq a$ and $a \in A(j)$ and $a' \in A(j')$ for $j,j' \geq 0$.  We shall prove that 
\[
|a-a'| \geq v_jg_j^{h-2} - g_j.
\]
Suppose not.  
If $j = j'$, then there exist $c,c' \in C(v_j,L)$ with $c \neq c'$ such that $a = p_j+c$ and $a' = p_j + c'$. 
By Lemma~\ref{Cassels:lemma:IntervalBasis} (iii) we have  
$|a-a'| = |c-c'| \geq v_jg_j^{h-2} - g_j.$  
Thus, if  $|a-a'| < v_j g_j^{h-2} - g_j$, then $j \neq j'$.  

The sequences $\{p_j\}_{j=0}^{\infty}$ and $\{g_j\}_{j=0}^{\infty}$ are strictly increasing sequences of positive integers.  If $j < j'$, then
$v_jg_j^{h-3} < v_{j'}g_{j'}^{h-3}$ and so 
\[
v_j g_j^{h-2} - g_j = ( v_jg_j^{h-3} -1 ) g_j 
< (v_{j'}g_{j'}^{h-3} - 1)g_{j'} =  v_{j'} g_{j'}^{h-2} - g_{j'}.
\]
Thus, if $j < j'$ and $|a-a'| < v_j g_j^{h-2} - g_j$, then also   $|a-a'| < v_{j'} g_{j'}^{h-2} - g_{j'}$.  Therefore, without loss of generality, we can assume that $j' < j$.  

By Lemma~\ref{Cassels:lemma:IntervalBasis} (ii) we have $a \geq p_j+g_j^h$ and $a' < p_{j'} + (h+L)g_{j'}^h.$
The inequality $|a-a'| < v_j g_j^{h-2} - g_j$ implies that 
\begin{align*}
a' & > a - v_jg_j^{h-2} + g_j > p_j + g_j^h - v_jg_j^{h-2} \\
&  = p_{j-1} + g_j^h  = p_{j-1} + 2^hg_{j-1}^h > p_{j-1} + hg_{j-1}^h .
\end{align*}
Combining the upper bound in Lemma~\ref{Cassels:lemma:IntervalBasis} (ii) with Lemma~\ref{Cassels:lemma:pjuj}, we get 
\[
a' < p_{j'} + (h+L)g_{j'}^h < (h+1+L)g_{j'}^h = 2^{2h}g_{j'}^h = 2^{h}g_{{j'}+1}^h = g_{{j'}+2}^h.
\]
Since $g_j^h < a' < g_{j' + 2}^h$, we see that $j' < j < j'+2$ and 
so $j = j' + 1$ and $a' = p_{j-1} + c'$ for some $c' \in C(v_{j-1},L)$ with $c' \geq hg_{j-1}^h$.  
By Lemma~\ref{Cassels:lemma:IntervalBasis} (ii),
we have $c' \equiv 0 \pmod{g_{j-1}^{h-1}}$ and so
\[
a' = p_{j-1} + c' \equiv p_{j-1} = p_j - v_j g_j^{h-2} \pmod{g_{j-1}^{h-1}}.
\]
Since
\[
g_{j-1}^{h-1} = 2^{h+j}g_{j-1}^{h-2} = 4v_j2^{h-2}g_{j-1}^{h-2} = 4v_jg_j^{h-2}
\]
it follows that 
\[
y = a' -  p_j \equiv - v_j g_j^{h-2} \pmod{ 4v_jg_j^{h-2} }.
\]
There exists $c \in C(v_j,L)$ such that $a = p_j + c$.  
Lemma~\ref{Cassels:lemma:IntervalBasis} (iv) implies that 
\[
|a-a'|= |c-(a'-p_j)| = |c-y|  \geq v_j g_j^{h-2} - g_j
\]
which is a contradiction.
This proves that if $a,a' \in A\setminus A(-1)$ with $a\neq a'$ and $a\in A(j)$, then $|a-a'| \geq v_jg_j^{h-2} - g_j.$

From Lemmas~\ref{Cassels:lemma:IntervalBasis} (ii) and~\ref{Cassels:lemma:pjuj} we also have  
\[
a = p_j + c < g_j^h + (h+L)g_j^h  = 2^{2h}g_j^h = (4g_j)^h
\]
and so $a^{(h-1)/h} < (4g_j)^{h-1}$
and
\[
\frac{|a-a'|}{a^{(h-1)/h}} > \frac{v_jg_j^{h-2}-g_j}{(4g_j)^{h-1}}
=  \frac{v_j}{4^{h-1}g_j}  - \frac{1}{4^{h-1}g_j^{h-2}}
=  \frac{1}{2^{3h-1}}  - \frac{1}{4^{h-1}g_j^{h-2}}.  
\]
Writing  $A$ as a strictly increasing sequence $A = \{a_k\}_{k=1}^{\infty}$ of nonnegative integers, we obtain 
\begin{align*}
\liminf_{k\rightarrow \infty}  \frac{a_{k+1}-a_k}{a_k^{(h-1)/h}} & \geq \liminf_{\substack{a,a' \in A\setminus A(-1) \\ a \neq  a'}} \frac{|a-a'|}{a^{(h-1)/h}} \\
& \geq  \liminf_{j\rightarrow \infty} 
\left( \frac{1}{2^{3h-1}} - \frac{1}{4^{h-1}g_j^{h-2}} \right) \\
& = \frac{1}{2^{3h-1}}.
\end{align*}
This completes the proof.  
\end{proof}

\bt[Cassels]      \label{Cassels:theorem:main3b}
For every integer $h \geq 3$  there exist a basis $C = \{c_n\}_{n=0}^{\infty}$ of order $h$ and real number $\lambda > 0$  such that $c_n = \lambda n^h + O\left(n^{h-1}\right)$.
\et

\begin{proof}
This follows immediately from Theorems~\ref{Cassels:theorem:main3a} and~\ref{Cassels:theorem:casimbed}.
\end{proof}

\boprob
Let $h \geq 2$.  Does there exist a basis $C = \{c_n\}_{n=0}^{\infty}$ of order $h$ such that $c_n = \gamma n^h + o(n^{h-1})$ for some $\gamma > 0$?
\eoprob

\boprob
Let $h \geq 2$.  Does there exist a basis $C = \{c_n\}_{n=0}^{\infty}$ of order $h$ such that $c_n = \gamma n^h + O(n^{h-2})$ for some $\gamma > 0$?
\eoprob

\boprob
Let $h \geq 2$.  
Compute or estimate 
\[
\sup \{ \lambda > 0 : \text{ there exists a basis $C = \{c_n\}_{n=0}^{\infty}$ of order $h$ such that $c_n \sim \lambda n^h$} \}.
\]
\eoprob

\section{Notes}
Raikov~\cite{raik37} and St\" ohr~\cite{stoh37}
independently constructed the first  examples of thin bases of order $h$.  
Another early, almost forgotten construction of thin bases is due to Chartrovsky~\cite{char40}.
The $g$-adic generalization of the Raikov-St\" ohr construction appears in work of Jia and Nathanson~\cite{jia96,nath89b} on minimal asymptotic bases.
The currently "thinest" bases of finite order appear in recent papers by Hofmeister~\cite{hofm01} and Blomer~\cite{blom03}.
An old but still valuable survey of combinatorial problems in additive number theory is St\" ohr~\cite{stoh55}.

The classical bases in additive number theory are the squares, cubes, and, for every integer $k \geq 4$, the $k$th powers of nonnegative integers, and also the sets of polygonal numbers and of prime numbers.  Using probability arguments, one can prove that all of the classical bases contain thin subsets that are bases of order $h$ for suffficiently large $h$ (Choi-Erd\H os-Nathanson~\cite{choi-erdo-nath80}, Erd\H os-Nathanson~\cite{nath81a}, Nathanson~\cite{nath81b}, Wirsing\cite{wirs86}, and Vu~\cite{vu00}).

The construction in this paper of polynomially asymptotic thin bases of order $h$ appeared in the classic paper of Cassels~\cite{cass57} in 1957.  
There is a recent quantitative improvement by Schmitt~\cite{schm06}, and also related work on  Cassels bases by Grekos, Haddad, Helou, and Pihko~\cite{grek-hadd-helo-pihk06} and Nathanson~\cite{nath09c}.

\def\cprime{$'$} \def\cprime{$'$} \def\cprime{$'$}
\providecommand{\bysame}{\leavevmode\hbox to3em{\hrulefill}\thinspace}
\providecommand{\MR}{\relax\ifhmode\unskip\space\fi MR }
\providecommand{\MRhref}[2]{%
  \href{http://www.ams.org/mathscinet-getitem?mr=#1}{#2}
}
\providecommand{\href}[2]{#2}

\end{document}